\magnification=\magstep1
\input amstex
\documentstyle{amsppt}
\voffset=-3pc
\baselineskip=12pt
\parskip=6pt
\def\Hom{\text{\bf Hom}}
\def\pHom{\text{\rm{Hom}}}
\def\GL{\text{\rm{GL}}}
\def\plain{\text{plain}}
\def\Ker{\text{\rm{Ker }}}
\def\Im{\text{\rm{Im }}}
\def\id{\text{id}}
\def\cO{\Cal O}
\def\cA{\Cal A}
\def\cB{\Cal B}
\def\cE{\Cal E}
\def\cF{\Cal F}
\def\cG{\Cal G}
\def\cC{\Cal C}
\def\cK{\Cal K}
\def\bC{\Bbb C}
\def\cR{\Cal R}

\def\cJ{\Cal J}
\def\fm{\frak m}

\topmatter
\title Coherent Sheaves and Cohesive Sheaves\endtitle\footnote""{\ Research 
partially supported by NSF grant
DMS0700281, the Mittag--Leffler Institute, and the Clay Institute.}
\author L\'aszl\'o Lempert\endauthor
\address Department of Mathematics, Purdue University, West Lafayette,
IN 47907\endaddress
\abstract
We consider coherent and cohesive sheaves of $\cO$--modules over open sets
$\Omega\subset\bC^n$.
We prove that coherent sheaves, and certain other sheaves derived from them, 
are cohesive; and conversely, certain sheaves derived from
cohesive sheaves are coherent.
An important tool in all this, also proved here, is that the sheaf of Banach
space valued holomorphic germs is flat.
\endabstract
\endtopmatter

\document
\noindent{\sl To Linda Rothschild on her birthday}
\head 1.\ Introduction\endhead

The theory of coherent sheaves has been central to algebraic and analytic geometry in the past fifty
years.
By contrast, in infinite dimensional analytic geometry coherence is irrelevant, as most sheaves
associated with infinite dimensional complex manifolds are not even finitely generated over the
structure sheaf, let alone coherent.
In a recent paper with Patyi, [LP], we introduced the class of so called cohesive sheaves in Banach
spaces, that seems to be the correct replacement of coherent sheaves---we were 
certainly able to show
that many sheaves that occur in the subject are cohesive, and for cohesive 
sheaves Cartan's Theorems A
and B hold.
We will go over the definition of cohesive sheaves in Section 2, but for a precise formulation of the
results above the reader is advised to consult [LP].

While cohesive sheaves were designed to deal with infinite dimensional 
problems, they make sense in
finite dimensional spaces as well, and there are reasons to study them in this 
context, too.
First, some natural sheaves even over finite dimensional manifolds are not 
finitely generated:\ for
example the sheaf $\cO^E$ of germs of holomorphic functions taking values in a fixed infinite
dimensional Banach space $E$ is not.
It is not quasicoherent, either (for this notion, see [Ha]), but it is cohesive.
Second, a natural approach to study cohesive sheaves in infinite dimensional manifolds would be to
restrict them to various finite dimensional submanifolds.

The issue to be addressed in this paper is the relationship between coherence and cohesion in finite
dimensional spaces.
Our main results are Theorems 4.3, 4.4, and 4.1.
Loosely speaking, the first says that coherent sheaves are cohesive, and the second that they remain
cohesive even after tensoring with the sheaf $\cO^F$ of holomorphic germs valued in a Banach space $F$.
A key element of the proof is that $\cO^F$ is flat, Theorem 4.1.
This latter is also relevant for the study of subvarieties.
On the other hand, Masagutov showed that $\cO^F$ is not free in general, see [Ms, Corollary 1.4].

The results above suggest two problems, whose resolution has eluded us.
First, is the tensor product of a coherent sheaf with a cohesive sheaf itself 
cohesive?
Of course, one can also ask the more ambitious question whether the tensor product of two cohesive
sheaves is cohesive, but here one should definitely consider some kind of ``completed'' tensor product,
and it is part of the problem to find which one.
The second problem is whether a finitely generated cohesive sheaf is coherent.
If so, then coherent sheaves could be defined as cohesive sheaves of finite type.
We could only solve some related problems:\ according to 
Corollary 4.2, any finitely generated
subsheaf of $\cO^F$ is coherent; and cohesive subsheaves of coherent sheaves are also
coherent, Theorem 5.4.

\head 2.\ Cohesive sheaves, an overview\endhead

In this Section we will review notions and theorems related to the theory of cohesive sheaves, following
[LP].
We assume the reader is familiar with very basic sheaf theory.
One good reference to what we need here---and much more---is [S].
Let $\Omega\subset\bC^n$ be an open set and $E$ a complex Banach space.
A function $f\colon\Omega\to E$ is holomorphic if for each $a\in\Omega$ there is a linear map
$L\colon\bC^n\to E$ such that
$$
f(z)=f(a)+L(z-a)+o \|z-a\|,\qquad z\to a.
$$
This is equivalent to requiring that in each ball $B\subset\Omega$ centered at any $a\in\Omega$ our $f$
can be represented as a locally uniformly convergent power series $f(z)=\sum_j c_j (z-a)^j$, with
$j=(j_1,\ldots,j_n)$ a nonnegative multiindex and $c_j\in E$.
We denote by $\cO^E_\Omega$ or just $\cO^E$ the sheaf of holomorphic $E$--valued germs over $\Omega$.
In particular, $\cO=\cO^\bC$ is a sheaf of rings, and $\cO^E$ is a sheaf of $\cO$--modules.
Typically, instead of a sheaf of $\cO$--modules we will just talk about $\cO$--modules.

\definition{Definition 2.1}The sheaves $\cO^E=\cO^E_\Omega\to\Omega$ are called plain sheaves.
\enddefinition

\proclaim{Theorem 2.2} [Bi, Theorem 4], [Bu, p.~331] or [L, Theorem 2.3].
If $\Omega\subset\bC^n$ is
pseudoconvex and $q=1,2,\ldots$, then $H^q(\Omega,\cO^E)=0$.
\endproclaim

Given another Banach space $F$, we write $\pHom(E,F)$ for the Banach space of continuous linear maps
$E\to F$.
If $U\subset\Omega$ is open, then any holomorphic $\Phi\colon U\to\pHom(E,F)$ induces a homomorphism
$\varphi\colon\cO^E|U\to\cO^F|U$, by associating with the germ of a holomorphic $e\colon V\to E$ at
$\zeta\in V\subset U$ the germ of the function $z\mapsto\Phi(z) e(z)$, again at $\zeta$.
Such homomorphisms and their germs are called plain.
The sheaf of plain homomorphisms between $\cO^E$ and $\cO^F$ is denoted $\Hom_{\plain}(\cO^E,\cO^F)$.
If $\Hom_{\cO}(\cA,\cB)$ denotes the sheaf of $\cO$--homomorphisms between $\cO$--modules $\cA$ and $\cB$,
then
$$
\Hom_{\plain}(\cO^E,\cO^F)\subset\Hom_{\cO}(\cO^E,\cO^F)\tag2.1
$$
is an $\cO$--submodule.
In fact, Masagutov showed that the two sides in (2.1) are equal unless $n=0$,
see [Ms, Theorem 1.1], but for the moment we do not need this.
The $\cO(U)$--module of sections $\Gamma(U,\Hom_{\plain}(\cO^E,\cO^F))$ is 
in one--to--one correspondence with the
$\cO(U)$--module $\pHom_{\plain}(\cO^E|U,\cO^F|U)$ of plain homomorphisms.
Further, any germ $\Phi\in\cO_z^{\pHom(E,F)}$ induces a germ $\varphi\in\Hom_{\plain}(\cO^E,\cO^F)_z$.
As pointed out in [LP, Section 2], the resulting map is an isomorphism
$$
\cO^{\pHom(E,F)}\overset\approx\to\rightarrow \Hom_{\plain}(\cO^E,\cO^F)\tag2.2
$$
of $\cO$--modules.

\definition{Definition 2.3}An analytic structure on an $\cO$--module $\cA$ is the choice, for each plain
sheaf $\cE$, of a submodule $\Hom(\cE,\cA)\subset\Hom_{\cO}(\cE,\cA)$, subject to
\item{(i)}if $\cE,\cF$ are plain sheaves and $\varphi\in\Hom_{\plain}(\cE,\cF)_z$ for some $z\in\Omega$,
then $\varphi^*\Hom(\cF,\cA)_z\subset\Hom(\cE,\cA)_z$; and
\item{(ii)}$\Hom(\cO,\cA)=\Hom_{\cO}(\cO,\cA)$.
\enddefinition

If $\cA$ is endowed with an analytic structure, one says that $\cA$ is an analytic sheaf.
The reader will realize that this is different from the traditional terminology, where ``analytic
sheaves'' and ``$\cO$--modules'' mean one and the same thing.

For example, one can endow a plain sheaf $\cG$ with an analytic structure by setting
$$
\Hom(\cE,\cG)=\Hom_{\plain}(\cE,\cG).
$$
Unless stated otherwise, we will always consider plain sheaves endowed with this analytic
structure.---Any $\cO$--module $\cA$ has two extremal analytic structures.
The maximal one is given by $\Hom(\cE,\cA)=\Hom_{\cO} (\cE,\cA)$.
In the minimal structure, $\Hom_{\min}(\cE,\cA)$ consists of germs $\alpha$ than can be written
$\alpha=\sum\beta_j\gamma_j$ with 
$$
\gamma_j\in\Hom_{\plain}(\cE,\cO)\quad\text{and}\quad\beta_j\in
\Hom_{\cO}(\cO,\cA),\quad j=1,\ldots,k.
$$

An $\cO$--homomorphism $\varphi\colon\cA\to\cB$ of $\cO$--modules
induces a homomorphism
$$
\varphi_*\colon\Hom_{\cO} (\cE,\cA)\to\Hom_{\cO}(\cE,\cB)
$$
for $\cE$ plain.
When $\cA$, $\cB$ are analytic sheaves, we say that $\varphi$ is analytic if
$$
\varphi_*
\Hom(\cE,\cA)
\subset\Hom(\cE,\cB)
$$ 
for all plain sheaves $\cE$.
It is straightforward to check that if $\cA$ and $\cB$ themselves are plain sheaves, then $\varphi$ is
analytic precisely when it is plain.
We write $\pHom(\cA,\cB)$ for the $\cO(\Omega)$--module of analytic homomorphisms $\cA\to\cB$ and $\Hom(\cA,\cB)$ for the sheaf of germs of analytic homomorphisms $\cA|U\to\cB|U$, with $U\subset\Omega$ open.
Again, one easily checks that, when $\cA=\cE$ is plain, this new notation is consistent with the one already in use.
Further,
$$
\pHom(\cA,\cB)\approx\Gamma(\Omega,\Hom(\cA,\cB)).\tag2.3
$$

\definition{Definition 2.4} Given an $\cO$-homomorphism
$\varphi:\cA\to\cB$ of $\cO$-modules,
any analytic structure on $\cB$ induces one on $\cA$ by the formula
$$
\Hom(\cE,\cA)=\varphi_*^{-1}\Hom(\cE,\cB).
$$
If $\varphi$ is an epimorphism, then any analytic structure on $\cA$ induces one on $\cB$ by the formula
$$
\Hom(\cE,\cB)=\varphi_*\Hom(\cE,\cA).
$$
\enddefinition
[LP, 3.4] explains this construction in the cases when $\varphi$ is the inclusion of a submodule
$\cA\subset\cB$ and when $\varphi$ is the projection on a quotient $\cB=\cA/\cC$.

Given a family $\cA_i$, $i\in I$, of analytic sheaves, an analytic structure is induced on the sum $\cA=\bigoplus\cA_i$.
For any plain $\cE$ there is a natural homomorphism
$$
\bigoplus_i\Hom_{\cO} (\cE,\cA_i)\to\Hom_{\cO} (\cE,\cA),
$$
and we define the analytic structure on $\cA$ by letting $\Hom(\cE,\cA)$ be the image of $\bigoplus\Hom(\cE,\cA_i)$.
With this definition, the inclusion maps $\cA_i\to\cA$ and the projections
$\cA\to\cA_i$ are analytic.

\definition{Definition 2.5}A sequence $\cA\to\cB\to\cC$ of analytic sheaves and homomorphisms over
$\Omega$ is said to be completely exact if for every plain sheaf $\cE$ and
every pseudoconvex
$U\subset\Omega$ the induced sequence 
$$
\pHom(\cE|U,\cA|U)\to\pHom(\cE|U,\cB|U)\to\pHom(\cE|U,\cC|U)
$$
is exact.
A general sequence of analytic homomorphisms is completely exact if every three--term subsequence is
completely exact.
\enddefinition

\definition{Definition 2.6}An infinite completely exact sequence
$$
\ldots\to\cF_2\to\cF_1\to\cA\to 0\tag2.4
$$
of analytic homomorphisms is called a complete resolution of $\cA$ if each $\cF_j$ is plain.
\enddefinition

When $\Omega$ is finite dimensional, as in this paper, complete resolutions can be defined more simply:

\proclaim{Theorem 2.7}Let
$$
\ldots\to\cF_2\overset{\varphi_2}\to\longrightarrow \ \cF_1
\overset{\varphi_1}\to\longrightarrow \ \cA\overset{\varphi_0}\to
\longrightarrow 0\tag2.5
$$
be an infinite sequence of analytic homomorphisms over $\Omega\subset\bC^n$,
with each $\cF_j$ plain.
If for each plain $\cE$ over $\Omega$ the induced sequence
$$
\ldots\to\Hom(\cE,\cF_2)\to\Hom(\cE,\cF_1)\to\Hom(\cE,\cA)\to 0\tag2.6
$$
is exact, then (2.5) is completely exact.
\endproclaim

\demo{Proof}Setting $\cE=\cO$ in (2.6) we see that (2.5) is exact.
Let $\cK_j=\Ker\varphi_j=\Im\varphi_{j+1}$, and endow it with the analytic structure induced by the
embedding $\cK_j\hookrightarrow\cF_j$, as in Definition 2.4.
The exact sequence 
$\cO\to\cK_j\hookrightarrow\cF_j\overset{\varphi_j}\to\longrightarrow\
\cK_{j-1}\to 0$ induces a sequence
$$
0\to\Hom(\cE,\cK_j)\to\Hom(\cE,\cF_j)\to\Hom(\cE,\cK_{j-1})\to 0,\tag2.7
$$
also exact since (2.6) was.
Let $U\subset\Omega$ be pseudoconvex.
Then in the long exact sequence associated with (2.7)
$$
\multline
\ldots\to H^q(U,\Hom(\cE,\cF_j))\to H^q (U,\Hom(\cE,\cK_{j-1}))\to \\
\to H^{q+1} (U,\Hom(\cE,\cK_j))\to H^{q+1} (U,\Hom(\cE,\cF_j))\to\ldots
\endmultline\tag2.8
$$
the first and last terms indicated vanish for $q\geq 1$ by virtue of 
Theorem 2.2 and (2.2).
Hence the middle terms are isomorphic:
$$
\multline
H^q(U,\Hom(\cE,\cK_{j-1}))\approx H^{q+1} (U,\Hom(\cE,\cK_j))\approx\ldots\\
\ldots\approx H^{q+n}(U,\Hom(\cE,\cK_{j+n-1}))\approx 0.
\endmultline
$$
Using this and (2.3), the first few terms of the
sequence (2.8) are
$$
0\to\pHom(\cE|U,\cK_j|U)\to\pHom(\cE|U,\cF_j|U)\to\pHom(\cE|U,\cK_{j-1}|U)
\to 0.
$$
The exactness of this latter implies 
$\ldots\to\pHom(\cE|U,\cF_1|U)\to\pHom(\cE|U,\cA|U)\to 0$ is exact, and so
(2.5) is indeed completely exact.
\enddemo

\definition{Definition 2.8}An analytic sheaf $\cA$ over $\Omega\subset\bC^n$ is cohesive if each $z\in\Omega$ has a neighborhood over which $\cA$ has a complete resolution.
\enddefinition

The simplest examples of cohesive sheaves are the plain sheaves, that have complete resolutions of form $\ldots\to 0\to 0\to\cE\to\cE\to 0$.
The main result of [LP] is the following generalization of Cartan's Theorems A and B, see Theorem 2 of the Introduction there:

\proclaim{Theorem 2.9}Let $\cA$ be a cohesive sheaf over a pseudoconvex $\Omega\subset\bC^n$.
Then
\itemitem{(a)}$\cA$ has a complete resolution over all of $\Omega$;
\itemitem{(b)}$H^q (\Omega,\cA)=0$ for $q\geq 1$. 
\endproclaim

\head 3.\ Tensor products\endhead

Let $R$ be a commutative ring with a unit and $A,B$ two $R$--modules.
Recall that the tensor product $A\otimes_R B=A\otimes B$ is the $R$--module 
freely generated by the set $A\times B$, modulo the submodule generated by 
elements of form
$$
(ra+a',b)-r(a,b)-(a',b)\quad\text{and}\quad
(a,rb+b')-r(a,b)-(a,b'),
$$
where $r\in R,\ a, a'\in A$, and $b, b'\in B$.
The class of $(a,b)\in A\times B$ in $A\otimes B$ is denoted $a\otimes b$.
Given homomorphisms $\alpha\colon A\to A'$, $\beta\colon B\to B'$ of $R$--modules, $\alpha\otimes\beta\colon A\otimes B\to A'\otimes B'$ denotes the unique homomorphism satisfying $(\alpha\otimes\beta)(a\otimes b)=\alpha(a)\otimes\beta(b)$.

A special case is the tensor product of Banach spaces $A, B$; here $R=\bC$.
The tensor product $A\otimes B$ is just a vector space, on which in general 
there are several natural ways to introduce a norm.
However, when $\dim A=k<\infty$, all those norms are equivalent, and turn $A\otimes B$ into a Banach space.
For example, if a basis $a_1,\ldots,a_k$ of $A$ is fixed, any $v\in A\otimes B$ can be uniquely written $v=\sum a_j\otimes b_j$, with $b_j\in B$.
Then $A\otimes B$ with the norm
$$
\|v\|=\max_j \|b_j\|_B
$$
is isomorphic to $B^{\oplus k}$.

Similarly, if $\cR$ is a sheaf of commutative unital rings over a topological 
space $\Omega$ and $\cA,\cB$ are $\cR$--modules, then the tensor product sheaf 
$\cA\otimes_\cR\cB=\cA\otimes\cB$ can be defined, see e.g.~[S].
The tensor product is itself a sheaf of $\cR$--modules, its stalks 
$(\cA\otimes\cB)_x$ are just the tensor products of 
$\cA_x$ and $\cB_x$ over $\cR_x$. Fix now an open $\Omega\subset\bC^n$, 
an $\cO$--module $\cA$, and an analytic sheaf $\cB$ over $\Omega$.
An analytic structure can be defined on $\cA\otimes\cB$ as follows.
For any plain sheaf $\cE$ there is a tautological $\cO$--homomorphism
$$
T=T_{\cE}\colon\cA\otimes\Hom(\cE,\cB)\to\Hom_{\cO} (\cE,\cA\otimes\cB),\tag3.1
$$
obtained by associating with $a\in\cA_\zeta$, $\epsilon\in\Hom(\cE,\cB)_\zeta$ first a section $\tilde a$ of $\cA$ over a neighborhood $U$ of $\zeta$, such that $\tilde a(\zeta)=a$; then defining $\tau^a\in\Hom_{\cO}(\cB,\cA\otimes\cB)_\zeta$ as the germ of the homomorphism
$$
\cB_z\ni b\mapsto\tilde a(z)\otimes b\in\cA_z\otimes\cB_z,\qquad z\in U;
$$
and finally letting $T(a\otimes\epsilon)=\tau^a\epsilon$.

\definition{Definition 3.1}The (tensor product) analytic structure on $\cA\otimes\cB$ is given by $\Hom(\cE,\cA\otimes\cB)=\Im T_\cE$.
\enddefinition

One quickly checks that this prescription indeed satisfies the axioms of an analytic structure. Equivalently, one can define
$\Hom(\cE,\cA\otimes\cB)\subset\Hom_\cO(\cE,\cA\otimes\cB)$
as the submodule spanned by germs of homomorphisms of the
form
$$
\cE|U @>\approx>>\cO\otimes\cE|U @>\alpha\otimes\beta>>
\cA\otimes\cB|U,
$$
where $U\subset\Omega$ is open, $\alpha:\cO|U\to\cA|U$ and
$\beta:\cE|U\to\cB|U$ are $\cO-$, resp. analytic homomorphisms
(and the first isomorphism is the canonical one).
The following is obvious.

\proclaim{Proposition 3.2}$T$ in (3.1) is natural:\ if $\alpha\colon\cA\to\cA'$ and $\beta\colon\cB\to\cB'$ are $\cO-$, resp.~analytic homomorphisms, then $T$ and the corresponding $T'$ fit in a commutative diagram
$$
\CD
\cA\otimes\Hom(\cE,\cB) @>{\alpha\otimes\beta_*}>> \cA'\otimes\Hom(\cE,\cB')\\
@VVTV @VVT'V\\
\Hom_{\cO} (\cE,\cA\otimes\cB) @>(\alpha\otimes\beta)_*>>\Hom_{\cO} (\cE,\cA'\otimes\cB').
\endCD
$$
\endproclaim

\proclaim{Corollary 3.3}If $\alpha,\beta$ are as above, then $\alpha\otimes\beta\colon\cA\otimes\cB\to\cA'\otimes\cB'$ is analytic.
\endproclaim

\proclaim{Proposition 3.4}If $\cA$ is an $\cO$--module and $\cB$ an analytic sheaf, then the tensor product analytic structure on $\cA\otimes\cO$ is the minimal one.
Further, the map
$$
\cB\ni b\mapsto 1\otimes b\in\cO\otimes\cB
$$
is an analytic isomorphism.
\endproclaim

Both statements follow from inspecting the definitions.

\proclaim{Proposition 3.5}If $\cA,\cA_i$ are $\cO$--modules and $\cB,\cB_i$ are analytic sheaves, then the obvious $\cO$--isomorphisms
$$
(\bigoplus\cA_i)\otimes\cB\overset\approx\to\longrightarrow \bigoplus
(\cA_i\otimes\cB),\quad
\cA\otimes (\bigoplus\cB_i)\overset\approx\to\longrightarrow\bigoplus
(\cA\otimes\cB_i)
$$
are in fact analytic isomorphisms.
\endproclaim

This follows from Definition 3.1, upon taking into account the distributive property of the tensor product of $\cO$--modules.
Consider now a finitely generated plain sheaf 
$\cF=\cO^F\approx\cO\oplus\ldots\oplus\cO$, with $\dim F=k$.
By putting together Propositions 3.4 and 3.5 we obtain analytic isomorphisms
$$
\cF\otimes\cB\approx (\cO\otimes\cB)\oplus\ldots\oplus
(\cO\otimes\cB)\approx\cB\oplus\ldots\oplus\cB.
$$
When $\cB=\cO^B$ is plain, this specializes to 
$$
\cO^F\otimes\cO^B\approx \cO^B\oplus\ldots\oplus\cO^B\approx 
\cO^{B^{\oplus k}}\approx\cO^{F\otimes B}. \tag3.2
$$

Later on we will need to know that inducing, in the sense of Definition 2.4, and tensoring are compatible.
Here we discuss the easy case, an immediate consequence of the tensor product being a right exact functor; the difficult case will have to wait until Section 6.

\proclaim{Proposition 3.6}Let $\psi\colon\cA\to\cA'$ be an epimorphism of $\cO$--modules and $\cB$ an analytic sheaf.
Then the tensor product analytic structure on $\cA'\otimes\cB$ is induced (in the sense of Definition 2.4) from the tensor product analytic structure on $\cA\otimes\cB$ by the epimorphism $\psi\otimes id_{\cB}\colon\cA\otimes\cB\to\cA'\otimes\cB$.
\endproclaim

\demo{Proof}We write $\cA\otimes\cB$, $\cA'\otimes\cB$ for the analytic sheaves endowed with the tensor product structure.
The claim means
$$
(\psi\otimes \id_{\cB})_*\Hom(\cE,\cA\otimes\cB)=\Hom(\cE,\cA'\otimes\cB)
$$
for every plain $\cE$.
But this follows from Definition 3.1 if we take into account the naturality of $T$ (Proposition 3.2) and that
$$
\psi\otimes \id_{\Hom(\cE,\cB)}\colon \cA\otimes\Hom(\cE,\cB)\to\cA'\otimes\Hom(\cE,\cB)
$$
is onto.
\enddemo

In the sequel it will be important to know when $T$ in (3.1) is injective.
This issue is somewhat subtle and depends on the analysis of Section 5.

\head 4.\ The main results\endhead

We fix an open set $\Omega\subset\bC^n$.
In the remainder of this paper all sheaves, unless otherwise stated, will be over $\Omega$.

\proclaim{Theorem 4.1}Let $\cF$ be a plain sheaf, $\cA\subset\cF$ finitely generated, and $\zeta\in\Omega$.
Then on some open $U\ni\zeta$ there is a finitely generated free subsheaf $\cE\subset\cF|U$ that contains $\cA|U$.
In particular, plain sheaves are flat.
\endproclaim

Recall that an $\cO$--module $\cF$ is flat if for every exact sequence $\cA\to\cB\to\cC$ of $\cO$--modules the induced sequence $\cA\otimes\cF\to\cB\otimes\cF\to\cC\otimes\cF$ is exact.

Theorem 4.1 and Oka's coherence theorem imply

\proclaim{Corollary 4.2}Finitely generated submodules of a plain sheaf are coherent.
\endproclaim

\proclaim{Theorem 4.3}A coherent sheaf, endowed with its minimal analytic structure, is cohesive.
\endproclaim

\proclaim{Theorem 4.4}If $\cA$ is a coherent sheaf and $\cB$ is a plain sheaf, then $\cA\otimes\cB$ is cohesive.
\endproclaim

Theorem 4.1 will be proved in Section 5, Theorems 4.3 and 4.4 in Section 7.

\head 5.\  Preparation\endhead

The main result of this Section is the following.
Throughout, $\Omega\subset\bC^n$ will be open.

\proclaim{Lemma 5.1}Let $P,Q$ be Banach spaces, $f\colon\Omega\to\pHom(P,Q)$ 
holomorphic, and $\zeta\in\Omega$.
\itemitem{(a)}If $\dim P<\infty$ then there are a finite dimensional 
$Q'\subset Q$, an open $U\ni\zeta$, and a holomorphic $q\colon U\to \GL(Q)$ 
such that $\Im q(z)f(z)\subset Q'$ for all $z\in U$.
\itemitem{(b)}If $\dim Q<\infty$ then there are a finite codimensional 
$P'\subset P$, an open $U\ni\zeta$, and a holomorphic $p\colon U\to \GL(P)$ 
such that $P'\subset\Ker f(z)p(z)$ for all $z\in U$.
\endproclaim

The proof depends on various extensions of the Weierstrass Preparation Theorem.
Let $A$ be a Banach algebra with unit $\bold 1$, and let $A^\times \subset A$ denote the open set of invertible elements.

\proclaim{Lemma 5.2}Let $f\colon\Omega\to A$ be holomorphic, $0\in\Omega$, and $d=0,1,2,\ldots$ such that 
$$
{\partial^j f\over\partial z_1^j}(0)=0\text\quad{ for}\quad j<d,\quad\text{ and }
\quad{\partial^d f\over\partial z_1^d} (0)\in A^\times.
$$
Then on some open $U\ni 0$ there is a holomorphic $\Phi\colon U\to A^\times$ such that, writing $z=(z_1,z')$
$$
\Phi(z)f(z)=\bold 1 z_1^d+\sum^{d-1}_{j=0} f_j (z') z_1^{j},\quad z\in U,
\tag5.1
$$
and $f_j(0)=0$.
\endproclaim

We refer the reader to [H\"o, 6.1].
The proof of Weierstrass' theorem given there for the case $A=\bC$ applies in this general setting as well.

\proclaim{Lemma 5.3}Let $0\in\Omega$, $E$ a Banach space, $E^*$ its dual,
$g\colon\Omega\to E$ (resp.~$h\colon\Omega\to E^*$) holomorphic functions 
such that
$$
{\partial^j g\over\partial z_1^j}\ (0)\neq 0\qquad (\text{resp. }{\partial^j h\over\partial z_1^j}\ (0)\neq 0),
\quad\text{for some }j.\tag5.2
$$
Then there are an open $U\ni 0$, a holomorphic $\Phi\colon U\to
\text{\rm GL}(E)$, and $0\neq e\in E$ (resp.~$0\neq e^* \in E^*$), such that
$$
\gathered
\Phi(z) g(z)=e z_1^d+\sum^{d-1}_{j=0} g_j (z') z_1^{j}\quad ,\quad z\in U,\\
(\text{resp. }h(z)\Phi(z)=e^* z_1^d+\sum^{d-1}_{j=0} h_j (z') z_1^{j}),
\endgathered\tag5.3
$$
with some $d=0,1,\ldots$, and $g_j(0)=0$, $h_j(0)=0$.
\endproclaim

\demo{Proof}We will only prove for $g$, the proof for $h$ is similar.
The smallest $j$ for which (5.2) holds will be denoted $d$.
Thus $\partial^d g/\partial z_1^d (0)=e\neq 0$.
Let $V\subset E$ be a closed subspace complementary to the line spanned by $e$, and define a holomorphic $f\colon\Omega\to\pHom (E,E)$ by
$$
f(z) (\lambda e+v)=\lambda g(z)+v z_1^d/d!,\quad \lambda\in\bC,\ v\in V.
$$
We apply Lemma 5.2 with the Banach algebra $A=\pHom (E,E)$; its invertibles form $A^\times=\text{GL}(E)$.
As
$$
{\partial^j f\over\partial z_1^j}\ (0)=0\text{ for }j<d\quad\text{and}
\quad{\partial^d f\over\partial z_1^d}\ (0)=\id_E,
$$
there are an open $U\ni 0$ and $\Phi\colon U\to \text{GL}(E)$ satisfying (5.1) and $f_j(0)=0$.
Hence
$$
\Phi(z) g(z)=\Phi(z) f(z)(e)=e z_1^d+\sum^{d-1}_{j=0}\ f_j(z')(e) z_1^{j},
$$
as claimed.
\enddemo

\demo{Proof of Lemma 5.1}We will only prove (a), part (b) is proved similarly.
The proof will be by induction on $n$, the case $n=0$ being trivial.

So assume the $(n-1)$--dimensional case and consider $\Omega\subset\bC^n$.
Without loss of generality we take $\zeta=0$.
Suppose first $\dim P=1$, say, $P=\bC$, and let $g(z)=f(z)(1)$.
Thus $g\colon\Omega\to Q$ is holomorphic.
When $g\equiv 0$ near 0, the claim is obvious; 
otherwise we can choose coordinates so that $\partial^j g/\partial z_1^j (0)\neq 0$ for some $j$.
By Lemma 5.3 there is a holomorphic $\Phi\colon U\to \text{GL}(Q)$ satisfying (5.3).
We can assume $U=U_1\times\Omega'\subset\bC\times\bC^{n-1}$.
Consider the holomorphic function $f'\colon\Omega'\to\pHom(\bC^{d+1},Q)$ given by
$$
f'(z')(\xi_0,\xi_1,\ldots,\xi_d)=e\xi_0+\sum_1^d g_j (z')\xi_j.
$$
By the inductive assumption, after shrinking $U$ and $\Omega'$, there are a $q'\colon U'\to \text{GL}(Q)$ and a finite dimensional $Q'\subset Q$ so that $\Im q'(z')f'(z')\subset Q'$ for all $z'\in U'$.
This implies $q'(z')\Phi(z) g(z)\in Q'$, and so with $q(z)=q'(z')\Phi(z)$ indeed $\Im q(z) f(z)\subset Q'$.

To prove the claim for $\dim P>1$ we use induction once more, this time on $\dim P$.
Assume the claim holds when $\dim P<k$, and consider a $k$--dimensional $P$, $k\geq 2$.
Decompose $P=P_1\oplus P_2$ with $\dim P_1=1$.
By what we have already proved, there are an open $U\ni 0$, a holomorphic $q_1\colon U\to \text{GL}(Q)$, and a finite dimensional $Q_1\subset Q$ such that $q_1(z)f(z)P_1\subset Q_1$.
Choose a closed complement $Q_2\subset Q$ to $Q_1$, and with the projection $\pi\colon Q_1\oplus Q_2\to Q_2$ let
$$
f_2 (z)=\pi q_1 (z) f(z)\in\pHom (P, Q_2).\tag5.4
$$
As $\dim P_2=k-1$, by the inductive hypothesis there are a finite dimensional $Q'_2\subset Q_2$ and (after shrinking $U$) a holomorphic $q'_2\colon U\to \text{GL}(Q_2)$ such that $q'_2(z)f_2(z)P_2\subset Q'_2$.
We extend $q'_2$ to $q_2\colon U\to \text{GL}(Q)$ by taking it to be the identity on $Q_1$.
Then $q_2(z)f_2(z)P_2\subset Q_2'$ and
$$
q_2 (z) q_1 (z) f(z) P_1\subset Q_1.\tag5.5
$$
Further, (5.4) implies $(q_1 (z) f(z)-f_2(z))P\subset Q_1$ and so
$$
q_2(z)q_1(z)f(z)P_2\subset q_2(z)Q_1+q_2(z)f_2(z)P_2\subset Q_1\oplus Q'_2.
\tag5.6
$$
(5.5) and (5.6) show that $q=q_2q_1$ and $Q'=Q_1\oplus Q'_2$ satisfy the requirements, and the proof is complete.
\enddemo

\demo{Proof of Theorem 4.1}Let $\cF=\cO^F$ and let $\cA$ be generated by holomorphic $f_1,\ldots,f_k\colon\Omega\to F$.
These functions define a holomorphic $f\colon\Omega\to\pHom(\bC^k,F)$ by
$$
f(z)(\xi_1,\ldots,\xi_k)=\sum_j\xi_j f_j(z).
$$
Choose a finite dimensional $Q'\subset F$, an open $U\ni\zeta$, and a holomorphic $q\colon U\to \text{GL}(F)$ as in Lemma 5.1(a).
Then $\cA'=q^{-1}\cO^{Q'}|U\subset\cF$ is finitely generated and free;
moreover, it contains the germs of each $f_j|U$, hence also $\cA|U$.

As to flatness:\ it is known, and easy, that the direct limit of flat modules is 
flat ([Mt, Appendix B]).
As each stalk of $\cF$ is the direct limit of its finitely generated free 
submodules, it is flat.
\enddemo

Here is another consequence of Lemma 5.1.

\proclaim{Theorem 5.4}Let $\cA$ be a coherent sheaf
and let $\cB\subset\cA$ be a submodule.
If there are a plain sheaf $\cO^F=\cF$ and an $\cO$--epimorphism $\varphi\colon \cF\to\cB$, then $\cB$ is coherent.
In particular, cohesive subsheaves of $\cA$ are coherent.
\endproclaim

If $\cF$ of the theorem is finitely generated, then so is $\cB$, and its coherence is immediate from the definitions.
For the proof of the general statement we need the notion of depth.
Recall that given an $\cO$--module $\cA$, the depth of a stalk 
$\cA_\zeta$ is 
$0$ if there is a submodule $0\neq L\subset\cA_\zeta$ 
annihilated by the maximal ideal $\fm_\zeta\subset\cO_\zeta$.
Otherwise depth $\cA_\zeta>0$.
(For the general notion of depth, see [Mt, p.~130]; the version we use here 
is the one e.g.~in [Ms, Proposition 4.2], at least in the
positive dimensional case.)

\proclaim{Lemma 5.5}If $\cA$ is a coherent sheaf, then 
$$
D=\{z\in\Omega\colon\text{ depth } \cA_z=0\}
$$
is a discrete set.
\endproclaim

\demo{Proof}Observe that, given a compact polydisc $K\subset\Omega$, 
the $\cO(K)$--module $\Gamma(K,\cA)$ is finitely generated.
Indeed, if $0\to\cA'|K\to\cO^{\oplus p}|K\to\cA|K\to 0$ is an exact sequence
of $\cO|K$--modules, then $H^1(K,\cA')=0$ implies that
$$
\cO(K)^{\oplus p}\approx\Gamma(K,\cO^{\oplus p})\to\Gamma (K,\cA)
$$
is surjective.
We shall also need the fact that $\cO(K)$ is Noetherian, see e.g.~[F].

As for the lemma, we can assume $\dim\Omega>0$.
If $z\in D$, there is a nonzero submodule $B\subset\cA_z$ such that 
$\frak m_z B=0$.
Let $\cB^z$ denote the skyscraper sheaf over $\Omega$ whose only nonzero 
stalk is $B$, at $z$.
We do this construction for every $z\in D$.
With $K\subset\Omega$ a compact polydisc, the submodule
$$
\sum_{z\in D\cap K}\Gamma(K,\cB^z)\subset\Gamma (K,\cA)\tag5.7
$$
is finitely generated.
But $\Gamma(K,\cB^z)\neq 0$ consists of (certain) sections of $\cA$ supported
at $z$.
It follows that the sum in (5.7) is a direct sum, hence in fact a finite 
direct sum.
In other words, $D\cap K$ is finite for every compact polydisc $K$, and $D$ 
must be discrete.
\enddemo

\demo{Proof of Theorem 5.4}We can suppose $\dim\Omega>0$.
First we assume that, in addition, depth $\cA_z>0$ for every 
$z$.
Since coherence
is a local property, and $\cA$ is locally finitely generated, we can assume that $\Omega$ is a ball, and there are a finitely 
generated plain sheaf $\cO^E=\cE\approx\cO\oplus\ldots\oplus\cO$ 
and an epimorphism $\epsilon\colon\cE\to\cA$.
According to Theorem 7.1 in [Ms], $\varphi$ factors through $\epsilon$:\ there is an $\cO$--homomorphism $\psi\colon\cF\to\cE$ such that $\varphi=\epsilon\psi$.
(Masagutov in his proof of Theorem 7.1 relies on a result of the present paper, but the reasoning is not circular.
What the proof of [Ms, Theorem 7.1] needs is our Theorem 4.3, whose proof is independent of Theorem 5.4 we are justifying here.)
Then [Ms, Theorem 1.1] implies $\psi$ is plain.

In view of Lemma 5.1(b), there are a finite codimensional $F'\subset F$ and
a plain isomorphism $\rho\colon\cF\to\cF$ such that $\psi\rho|\cO^{F'}=0$.
If $F''\subset F$ denotes a (finite dimensional) complement to $F'$, then 
$\psi\rho(\cO^{F''})=\psi\rho(\cF)=\psi(\cF)$.
Hence $\epsilon\psi\rho(\cO^{F''})=\epsilon\psi(\cF)=\cB$ is finitely generated ; as a submodule of a coherent sheaf, itself must be coherent.

Now take an $\cA$ whose depth is $0$ at some $z$.
In view of Lemma 5.5 we can assume that there is a single such $z$.
With
$$
C=\{a\in\cA_z\colon{\frak m}_z^k a=0\text{ for some }k=1,2,\ldots\},
$$
let $\cC\subset\cA$ be the skyscraper sheaf over $\Omega$ whose only nonzero stalk is $C$, at $z$.
As $C$ is finitely generated, $\cC$ and $\cA/\cC$ are coherent.
Also, depth$(\cA/\cC)_\zeta>0$ for every $\zeta\in\Omega$.
Therefore by the first part of the proof $\cB/\cB\cap\cC\subset\cA/\cC$ is 
coherent.
Since $\cB\cap\cC$, supported at the single point $z$, is coherent, the Three Lemma implies $\cB$ is coherent, as claimed.
\enddemo

\head 6.\ $\Hom$ and $\otimes$.\endhead
The main result of this Section is the following. Let $\cA$ be an
$\cO$--module and $\cB$ an analytic sheaf.
Recall that, given a plain sheaf $\cE$, in Section 3 we introduced a 
tautological $\cO$--homomorphism
$$
T=T_\cE\colon\cA\otimes\Hom(\cE,\cB)\to\Hom_{\cO}(\cE,\cA\otimes\cB),\tag6.1
$$
and $\Hom(\cE,\cA\otimes\cB)$ was defined as the image of $T_\cE$.

\proclaim{Theorem 6.1}If $\cB$ is plain, then $T$ is injective.
\endproclaim

Suppose $\cE=\cO^E,\cB=\cO^B$ are plain sheaves. If $\zeta\in\Omega$, a
$\bC$--linear map
$$
S^\zeta\colon\cA_\zeta\otimes\cO_\zeta^{\pHom(E,B)}\to
\pHom_\bC(E,\cA_\zeta\otimes\cO_\zeta^B)
\tag6.2
$$
can be defined as follows. Let $a\in\cA_\zeta$,
$\Theta\in\cO_\zeta^{\pHom(E,B)}$, then
$$
S^\zeta(a\otimes\Theta)(e)=a\otimes\Theta e,\qquad e\in E,
$$
where on the right $e$ is thought of as a constant germ $\in\cO_\zeta^E$.
The key to Theorem 6.1 is the following

\proclaim{Lemma 6.2}Let $E,B$ be Banach spaces, $\zeta\in\Omega$, and let
$M$ be an $\cO_\zeta$--module. Then the tautological homomorphism
$$
S\colon M\otimes\cO_\zeta^{\pHom(E,B)}\to
\pHom_\bC(E,M\otimes\cO_\zeta^B)\tag6.3
$$
given by $S(m\otimes\Theta)(e)=m\otimes\Theta e$, for $e\in E$, 
is injective.
\endproclaim

As $T$ was, $S$ is also natural with respect to $\cO_\zeta$--homomorphisms 
$M\to N$.
The claim of the lemma is obvious when $M$ is free, for then tensor products 
$M\otimes L$ are just direct sums of copies of $L$.
The claim is also obvious when $M$ is a direct summand in a free module 
$M'=M\oplus N$, as the tautological homomorphism for $M'$ decomposes into 
the direct sum of the 
tautological homomorphisms for $M$ and $N$.
 
\demo{Proof of Lemma 6.2}The proof is inspired by the proof of [Ms, Theorem 1.3].
The heart of the matter is to prove when $M$ is finitely generated.
Let us write $(L_n)$ for the statement of the lemma for $M$ finitely generated and
$n=\dim\Omega$; we prove it by induction on $n$.
$(L_0)$ is trivial, as $\cO_\zeta\approx\bC$ is a field and any module over it is free.
So assume $(L_{n-1})$ for some $n\geq 1$, and prove $(L_n)$.
We take $\zeta=0$.

Step $1^\circ$.
First we verify $(L_n)$ with the additional assumption that $gM=0$ with some 
$0\neq g\in\cO_0$.
By Weierstrass' preparation theorem we can take $g$ to be (the germ of) a 
Weierstrass polynomial of degree $d\geq 1$ in the $z_1$ variable.
We write $z=(z_1,z')\in\bC^n$, and $\cO'_0,\cO^{'F}_0$ for the ring/module of
the corresponding germs in $\bC^{n-1}$ (here $F$ is any Banach space).
We embed $\cO'_0\subset\cO_0$, $\cO^{'F}_0\subset\cO_0^F$ as germs independent 
of $z_1$.
This makes $\cO_0$--modules into $\cO'_0$--modules.
In the proof tensor products both over $\cO_0$ and $\cO'_0$ will occur; 
we keep writing $\otimes$ for the former and will write $\otimes'$ for the 
latter.

We claim that the $\cO'_0$--homomorphism
$$
i:M\otimes'\cO^{'F}_0\to M\otimes\cO_0^F,\qquad
i(m\otimes' f')=m\otimes f',
$$
is in fact an isomorphism. To verify it is surjective, consider
$m\otimes f\in M\otimes\cO^F_0$. 
By Weierstrass' division theorem, valid for vector valued functions
as well (e.g., the proof in [GuR, p.~70] carries over verbatim),
$f$ can be written
$$
f=f_0g+\sum^{d-1}_{j=0} f'_j z_1^j,
\qquad f_0\in\cO_0^F, f'_j\in\cO^{'F}_0.
$$
Thus $m\otimes f=m\otimes(f_0g+\sum f_j'z_1^j)=
i(\sum z_1^jm\otimes' f_j')$ is indeed in Im $i$. Further,
injectivity is clear if dim $F=k<\infty$, as
$M\otimes'\cO_0^{'F}\approx M^{\oplus k}$,
$M\otimes\cO_0^F\approx M^{\oplus k}$, and $i$ corresponds
to the identity of $M^{\oplus k}$.  For a general $F$ consider
a finitely generated submodule $A\subset\cO_0^{'F}$. Lemma 5.1(a)
implies that there are a neighborhood $U$ of $0\in\bC^{n-1}$,
a finite dimensional subspace $G\subset F$, and a holomorphic
$q:U\to\GL(F)$ such that the automorphism $\varphi'$ of 
$\cO_0^{'F}$ induced by $q$ maps $A$ into 
$\cO_0^{'G}\subset\cO_0^{'F}$. (The reasoning is the same
as in the proof of Theorem 4.1.) If extended to $\bC\times U$
independent of
$z_1$, $q$ also induces an automorphism $\varphi$ of $\cO_0^F$,
and $i$ intertwines the automorphisms id$_M\otimes'\varphi'$
and id$_M\otimes\varphi$. Now $i$ is injective between
$M\otimes'\cO_0^{'G}$ and $M\otimes\cO_0^G\subset M\otimes\cO_0^F$,
because dim $G<\infty$. It follows that $i$ is also injective
on the image of $M\otimes' A$ in $M\otimes'\cO^{'F}_0$. Since the
finitely generated $A\subset\cO_0^{'F}$ was 
arbitrary, $i$ itself is injective.

Applying this with $F=\pHom(E,B)$ and $F=B$, we obtain a commutative diagram
$$
\CD
M\otimes'\cO_0^{'\pHom(E,B)} @>\approx>> 
M\otimes\cO_0^{\pHom(E,B)}
\\
@VS'VV @VSVV \\
\pHom_\bC(E,M\otimes'\cO_0^{'B}) @>\approx>>
\pHom_\bC(E, M\otimes\cO_0^{B}).
\endCD
$$
Here $S'$ is also a tautological homomorphism.
Now $M$ is finitely generated over $\cO_0/g\cO_0$, and this
latter is a finitely generated $\cO'_0$-algebra by Weierstrass
division. It follows that $M$ is finitely generated 
over $\cO'_0$; by the induction hypothesis $S'$ is injective, 
hence so must be $S$.

Step $2^\circ$.
Now take an arbitrary finitely generated $M$.
Let $\mu\colon L\to M$ be an epimorphism from a free
finitely generated $\cO_0$--module $L$, and $K=\Ker\mu$.
If the exact sequence 
$$
0\to K \overset\lambda\to\rightarrow\ L\overset \mu\to\rightarrow\ M\to 0
\tag6.4
$$
splits, then $M$ is a direct summand in $L$ and, as said, the claim is 
immediate.
The point of the reasoning to follow is that, even if (6.4) does not split, it
does split up to torsion in the following sense:\ there are 
$\sigma\in\pHom(L,K)$ and $0\neq g\in\cO_0$ such that 
$\sigma\lambda\colon K\to K$ is multiplication by $g$.
To see this, let $Q$ be the field of fractions of $\cO_0$, and note that the 
induced linear map $\lambda_Q\colon K\otimes Q\to L\otimes Q$ of $Q$--vector 
spaces has a left inverse $\tau$.
Clearing denominators in $\tau$ then yields the $\sigma$ needed.

\hfuzz=300pt
We denote the tautological homomorphisms (6.3) for $K,L,M$ by $S_K, S_L, S_M$.
Tensoring and $\pHom$--ing (6.4) gives rise to a commutative diagram
$$
\minCDarrowwidth{0.4cm}
\CD
K\otimes\cO_0^{\pHom(E,B)}
@>\lambda_t>\overset\longleftarrow\to{\sigma_t}>
L\otimes\cO_0^{\pHom(E,B)} @>\mu_t>> M\otimes\cO_0^{\pHom(E,B)} @>>>
0\\
@VVS_KV @VVS_LV @VVS_MV\\
\pHom_{\bC} (E, K\otimes\cO_0^B)
@>\lambda_h>{\overset\longleftarrow\to{\sigma_h}}>
\pHom_{\bC} (E,L\otimes\cO_0^B) @>\mu_h>> 
\pHom_{\bC} (E,M\otimes\cO_0^B)@>>>0  
\endCD
$$
with exact rows.
Here $\lambda_t,\lambda_h$, etc.~just indicate homomorphisms induced on 
various modules by $\lambda$, etc.
Consider an element of $\Ker S_M$; it is of form $\mu_t u$, 
$u\in L\otimes\cO_0^{\pHom(E,B)} $.
Then $S_L u\in\Ker\mu_h=\Im\lambda_h$.
Let $S_L u=\lambda_h v$.
We compute
$$
S_L\lambda_t\sigma_t u=\lambda_h S_K\sigma_t u=\lambda_h\sigma_h S_L u=
\lambda_h\sigma_h\lambda_h v=\lambda_hgv=S_L gu.
$$
Since $L$ is free, $S_L$ is injective, so $\lambda_t\sigma_t u=gu$ 
and $g\mu_t u=\mu_t\lambda_t\sigma_t u=0$.
We conclude that $g\Ker S_M=0$.
Let $N\subset M$ denote the submodule of elements annihilated by $g$ and, for 
brevity,
set $H=\cO_0^{\pHom(E,B)}$, a flat module.
Multiplication by $g$ is a monomorphism on $M/N$, so the same holds on 
$M/N\otimes H$.
The exact sequence
$$
0\to N\otimes H\hookrightarrow M\otimes H\to M/N\otimes H\to 0
$$
then shows that in $M\otimes H$ the kernel of multiplication by $g$ is 
$N\otimes H$.
Therefore $N\otimes H\supset\Ker S_M$, and $\Ker S_M\subset\Ker S_N$.
But $gN=0$, so from Step $1^\circ$ it follows that $\Ker S_N=0$, and 
again $\Ker S_M=0$.

Step $3^\circ$. Having proved the lemma for finitely generated modules, 
consider an arbitrary $\cO_0$--module $M$. The inclusion 
$\iota:N\hookrightarrow M$
of a finitely generated submodule induces a commutative diagram
$$\CD
N\otimes\cO_0^{\pHom(E,B)}@>\iota_t>>M\otimes\cO_0^{\pHom(E,B)} \\
@VS_NVV @VSVV \\
\pHom_\bC(E,N\otimes\cO_0^B)@>\iota_h>>\pHom_\bC(E,M\otimes\cO_0^B),
\endCD
$$
with $S_N$ the tautological homomorphism for $N$. Flatness implies that
$\iota_t$, $\iota_h$ are injective; as $S_N$ is also injective by
what we have proved so far, $S$ itself is injective on the range of
$\iota_t$. As $N$ varies, these ranges cover all of
$M\otimes\cO_0^{\pHom(E,B)}$, hence $S$ is indeed injective.
\enddemo

\demo{Proof of Theorem 6.1} For $\zeta\in\Omega$ we embed 
$E\to\cO_\zeta^E$ as constant germs; this induces a $\bC$--linear map
$$
\rho:\Hom_\cO(\cO^E,\cA\otimes\cO^B)_\zeta\to
\pHom_\bC(E,\cA_\zeta\otimes\cO^B_\zeta).
$$
It will suffice to show that if we restrict $T$ to the stalk at $\zeta$
and compose it with $\rho$, the resulting map
$$
T^\zeta:\cA_\zeta\otimes\Hom(\cO^E,\cO^B)_\zeta\to
\pHom_\bC(E,\cA_\zeta\otimes\cO^B_\zeta),
$$
given by $T^\zeta(a\otimes\theta)(e)=a\otimes\theta e$,
is injective. But, by the canonical isomorphism
$\cO_\zeta^{\pHom(E,B)}@>\approx>>\Hom(\cO^E,\cO^B)_\zeta$, cf. (2.2),
$T^\zeta$ is injective precisely when $S^\zeta$ of (6.2) is; so that
Lemma 6.2 finishes off the proof.
\enddemo

Now we can return to the question how compatible are inducing in the sense of Definition 2.4 and tensoring.

\proclaim{Lemma 6.3}If $0\to\cA'\overset\varphi\to\rightarrow 
\cA\overset\psi\to\rightarrow \cA''\to 0$ is an exact sequence of 
$\cO$--modules and $\cB$ is a plain sheaf, then 
$\varphi\otimes\id_\cB$, resp.~$\psi\otimes\id_\cB$, induce from 
$\cA\otimes\cB$ the tensor product analytic structure on $\cA'\otimes\cB$, 
resp.~$\cA''\otimes\cB$.
\endproclaim

\demo{Proof}The case of $\cA''\otimes\cB$, in greater generality, is the content of Proposition 3.6.
Consider $\cA'\otimes\cB$.
Meaning by $\cA'\otimes\cB$ etc.~the analytic sheaves endowed with the tensor product structure, in light of Definition 2.4 we are to prove 
$$
\Hom(\cE,\cA'\otimes\cB)=
(\varphi\otimes\id_\cB)_*^{-1}\Hom(\cE,\cA\otimes\cB)\tag6.5
$$
for every plain $\cE$.
Again using that $\cB$ and $\Hom(\cE,\cB)$ are flat, from 
$0\to \cA'\to\cA\to\cA''\to 0$ we obtain a commutative diagram with exact 
rows
$$
\minCDarrowwidth{.5cm}
\CD
0 @>>> \cA'\otimes\Hom(\cE,\cB) @>\varphi_t>> \cA\otimes\Hom(\cE,\cB) @>\psi_t>> \cA''\otimes\Hom (\cE,\cB)\\
& & @VV T'V @VV TV @VV T''V\\
0 @>>> \Hom_{\cO} (\cE,\cA'\otimes\cB) @>\varphi_h>> \Hom_{\cO} (\cE,\cA\otimes\cB) @>\psi_h>> \Hom_{\cO} (\cE,\cA''\otimes\cB).
\endCD
$$
The vertical arrows are the respective tautological homomorphisms, and 
$\varphi_t=\varphi\otimes\id_{\Hom(\cE,\cB)}$, 
$\varphi_h=(\varphi\otimes\id_\cB)_*$, etc.
denote maps induced by $\varphi$, etc.
From this diagram, the left hand side of (6.5), $\Im T'$, 
is clearly 
contained in $\varphi_h^{-1}\Im T$, i.e.~in the right hand side.
To show the converse, suppose $\epsilon\in\Hom_{\cO}(\cE,\cA'\otimes\cB)$ is 
in $\varphi_h^{-1}\Im T$, say, $\varphi_h\epsilon=Tu$ with 
$u\in\cA\otimes\Hom (\cE,\cB)$.
Then $T''\psi_t u=\psi_h Tu=\psi_h\varphi_h\epsilon=0$.
Since $T''$ is injective by Theorem 6.1, $\psi_t u=0$.
It follows that $u=\varphi_t v$ with some $v\in\cA'\otimes\Hom(\cE,\cB)$, whence $\varphi_h T' v=T\varphi_t v=Tu=\varphi_h\epsilon$.
As $\varphi_h$ is also injective, $\epsilon=T' v$; that is, $\varphi_h^{-1}\Im T\subset\Im T'$, as needed.
\enddemo

\heading 7.\ Coherence and cohesion\endheading

\demo{Proof of Theorems 4.3 and 4.4}We have to show that if $\cA$ is a coherent sheaf and $\cB=\cO^B$ plain then $\cA\otimes\cB$ is cohesive.
This would imply that $\cA\otimes\cO$ is cohesive, and in view of Proposition 3.4 that $\cA\approx\cA\otimes\cO$, with its minimal analytic structure, is also cohesive.

We can cover $\Omega$ with open sets over each of which $\cA$ has a resolution by finitely generated free $\cO$--modules.
We can assume that such a resolution
$$
\ldots\to\cF_2\to \cF_1\to\cA\to 0
$$
exists over all of $\Omega$, and $\cF_j=\cO^{F_j}$, dim $F_j<\infty$.
If $\cE=\cO^E$ is plain then
$$
\ldots\rightarrow \cF_2\otimes\Hom(\cE,\cB) \rightarrow 
\cF_1\otimes\Hom(\cE,\cB) \rightarrow \cA\otimes\Hom(\cE,\cB)@>>> 0
$$
is also exact, $\Hom(\cE,\cB)\approx\cO^{\pHom(E,B)}$ being flat.
By Theorem 6.1 this sequence is isomorphic to
$$
\ldots\to\Hom (\cE,\cF_2\otimes\cB)\to\Hom (\cE,\cF_1\otimes\cB)\to\Hom(\cE,\cA\otimes\cB)\to 0,
$$
which then must be exact. Here 
$\cF_j\otimes\cB\approx \cO^{F_j\otimes B}$ analytically, cf. (3.2).
Now Theorem 2.7 applies. We conclude that
$$
\ldots @>>>\cF_2\otimes\cB\to\cF_1\otimes\cB\to\cA\otimes\cB\to 0
$$
is completely exact, and $\cA\otimes\cB$ is indeed cohesive.
\enddemo

\heading 8.\ Application.
Complex analytic subspaces and subvarieties\endheading

The terminology in the subject indicated in the title is varied and 
occasionally ambiguous, even in finite dimensional complex geometry.
Here we will use the terms ``complex subspace'' and ``subvariety'' to mean
different things. Following [GrR], a
complex subspace $A$ of an open $\Omega\subset\bC^n$ is obtained from a 
coherent subsheaf $\cJ\subset\cO$.
The support $|A|$ of the sheaf $\cO/\cJ$, endowed with the sheaf of rings 
$(\cO/\cJ)\big| |A|=\cO_A$ defines a ringed space, and the pair $(|A|,\cO_A)$ 
is the complex subspace in question.

For infinite dimensional purposes this notion is definitely not adequate, and 
in the setting of Banach spaces in [LP] we introduced a new notion that we 
called subvariety.
Instead of coherent sheaves, they are defined in terms of cohesive sheaves, furthermore, one has to specify a subsheaf $\cJ^E\subset\cO^E$ for each Banach space $E$, (thought of as germs vanishing on the subvariety), not just one $\cJ\subset\cO$.
The reason this definition was made was to delineate a class of subsets in Banach spaces that arise in complex analytical questions, and can be studied using complex analysis.
At the same time, the definition makes sense in $\bC^n$ as well, and it is natural to ask how subvarieties and complex subspaces in $\bC^n$ are related.
Before answering we have to go over the definition of subvarieties, following [LP].

An ideal system over $\Omega\subset\bC^n$ is the specification, for every 
Banach space $E$, of a submodule $\cJ^E\subset\cO^E$, subject to the following:
given $z\in\Omega$, $\varphi\in\cO_z^{\pHom(E,F)}$, and $e\in\cJ_z^E$, 
we have $\varphi e\in\cJ_z^F$.

Within an ideal system the support of $\cO^E/\cJ^E$ is the same for every $E\neq (0)$, and we call this set the support of the ideal system.

A subvariety $S$ of $\Omega$ is given by an ideal system of cohesive subsheaves $\cJ^E\subset\cO^E$.
The support of the ideal system is called the support $|S|$ of the subvariety, and we endow it with the sheaves $\cO_S^E=\cO^E/\cJ^E\big| |S|$ of modules over $\cO_S=\cO_S^{\bC}$.
The ``functored space'' $(|S|,\ E\mapsto\cO_S^E)$ is the subvariety $S$ in question.

\proclaim{Theorem 8.1}There is a canonical way to associate a subvariety with a complex subspace of $\Omega$  and vice versa.
\endproclaim

\demo{Proof, or rather construction}Let $i\colon\cJ\hookrightarrow\cO$ be 
the inclusion of a coherent sheaf $\cJ$ that defines a complex subspace. The 
ideal system $\cJ^E=\cJ\cO^E\subset\cO^E$ then gives rise to a 
subvariety, provided $\cJ^E$ with the analytic structure 
inherited from $\cO^E$ is cohesive.
Consider the diagram
$$
\CD
\cJ\otimes\cO^E @>i\otimes\id_{\cO^E}>>\cO\otimes\cO^E\\
@V\mu VV @VV\approx V\\
\cJ^E @>>>\cO^E.
\endCD
$$
Here the vertical arrow on the right, given by $1\otimes e\mapsto e$, is an analytic isomorphism by Proposition 3.4.
The vertical arrow $\mu$ on the left is determined by the commutativity of the diagram; it is surjective.
As $\cO^E$ is flat, $i\otimes\id_{\cO^E}$ is injective, therefore $\mu$ is an isomorphism.
If $\cJ\otimes\cO^E$ is endowed with the analytic structure induced by $i\otimes\id_{\cO^E}$, $\mu$ becomes an analytic isomorphism.
On the other hand, the induced structure of $\cJ\otimes\cO^E$ agrees with the tensor product analytic structure by Lemma 6.3 (set $\cA'=\cJ,\cA=\cO$), hence it is cohesive by Theorem 4.4.
The upshot is that $\cJ^E$ is indeed cohesive.

As to the converse, suppose $\cJ^E$ is a cohesive ideal system defining a subvariety.
Then $\cJ=\cJ^\bC\subset\cO$ is coherent by Theorem 5.4, and gives rise to a complex subspace. 
\enddemo

Theorem 8.1 is clearly not the last word on the matter.
First, it should be decided whether the construction in the theorem is a bijection between subvarieties and complex subspaces; second, the functoriality properties of the construction should be investigated.


\Refs
\widestnumber\key{GuR}
\ref\key Bi\by E.~Bishop\paper Analytic functions with values in a Fr\'echet space\jour Pacific J.~Math.
\vol 12\yr 1962\pages 1177--1192\endref

\ref\key Bu\by
L.~Bungart\paper Holomorphic functions with values in locally convex spaces and applications to integral formulas\jour Trans.~Amer.~Math.~Soc\vol 111\yr 1964
\pages 317--344\endref

\ref\key F\by J. Frisch\paper Points de platitude d'un morphisme d'espaces
analytiques complexes\jour Invent. Math.\vol 4\yr 1967\pages 118--138\endref

\ref\key GrR\by
H.~Grauert, R.~Remmert\book Theory of Stein spaces\publ Springer\publaddr Berlin\yr1979\endref

\ref\key GuR\by
R.~Gunning, H.~Rossi\book Analytic functions of several complex variables
\publ Prentice Hall\publaddr Englewood Cliffs, N.J.\yr 1965\endref

\ref\key Ha\by R.~Harthshorne\book Algebraic Geometry\publ Springer\publaddr New York\yr1977\endref

\ref\key H\"o\by
L.~H\"ormander\book An introduction to complex analysis in several variables\publ North Holland\publaddr Amsterdam\yr1973\endref

\ref\key L\by
J.~Leiterer\paper Banach coherent analytic Fr\'echet sheaves\jour Math.~Nachr.\vol 85\yr 1978\pages 91--109\endref

\ref\key LP\by
L.~Lempert, I.~Patyi\paper Analytic sheaves in Banach spaces\jour Ann.~Sci.~\'Ecole Norm.~Sup.\vol40\yr2007\pages 453--486\endref

\ref\key Ms\by V.~Masagutov\paper Homomorphisms of infinitely generated sheaves, manuscript\endref

\ref\key Mt\manyby
H.~Matsumura\book Commutative ring theory\publ Cambridge University Press
\publaddr Cambridge\yr1986\endref

\ref\key S\by J.P.~Serre\paper Faisceaux alg\'ebriques coh\'erents
\jour Ann.~Math.\vol 61\yr 1955\pages 197--278\endref
\endRefs
\enddocument